\documentclass[a4paper,11pt]{article}

\usepackage{amssymb,amsthm,amscd}
\usepackage{amsmath}
\usepackage{graphics}
\usepackage{latexsym,epic,eepic,epsfig}
\usepackage{color}

\usepackage{tocloft}

\usepackage{enumitem}

\usepackage{apptools}

\setitemize{noitemsep,leftmargin=15pt, topsep=0pt,parsep=0pt,partopsep=10pt}
\setenumerate{noitemsep,leftmargin=25pt, topsep=0pt,parsep=0pt,partopsep=0pt}

\theoremstyle{plain}
\newtheorem{theorem}{Theorem}[section]
\newtheorem{defn}[theorem]{Definition}
\newtheorem{prop}[theorem]{Proposition}
\newtheorem{la}[theorem]{Lemma}
\newtheorem{cor}[theorem]{Corollary}
\newtheorem{rk}[theorem]{Remark}

\AtAppendix{\counterwithin{theorem}{section}}

\def\interior{{\rm int\ }}

\def\barr{\begin{array}}
\def\earr{\end{array}}
\def\beqarr*{\begin{eqnarray*}}
\def\eeqarr*{\end{eqnarray*}}

\def\mapright#1{\smash{\mathop{\longrightarrow}\limits^{#1}}}
\def\maprdown#1{\Big\downarrow
                   \rlap{$\vcenter{\hbox{$\scriptstyle#1$}}$}}
\def\mapldown#1{\llap{$\vcenter{\hbox{$\scriptstyle#1$}}$}\Big\downarrow}

\newcommand{\C}{\mathbb{C}}
\newcommand{\D}{\mathbb{D}}
\newcommand{\N}{\mathbb{N}}

\newcommand{\Z}{\mathbb{Z}}

\begin{document}

\pagenumbering{arabic}

%\makeatletter

\title{\Large\bf  Holomorphic Shadowing for H\'{e}non Maps Revisited: an Implicit Function Theorem Perspective}

\author{Yi-Chiuan Chen\\
   Institute of  Mathematics, Academia Sinica, Taipei 106319, Taiwan\\
   Email: YCChen@math.sinica.edu.tw \\
 }

%\date{}

\maketitle

\begin{abstract}
 In studying the complex H\'enon maps, Mummert (in ``Holomorphic shadowing for H\'{e}non maps"  {\it Nonlinearity} {\bf 21} pp. 2887-2898, 2008) defined an operator the fixed points of which give rise to bounded orbits. This enabled him to obtain an estimate of the solenoid locus. Instead of the contraction mapping theorem, in this paper we present an implicit function theorem version of his result by providing an alternative proof of a theorem of Hubbard and Oberste-Vorth (in {\it Real and Complex Dynamical Systems}, pp.89-132, 1995).
\end{abstract}

\vspace{2cm}


\vspace{1cm}


%\makeatother

% \tableofcontents

\newpage

\section{Introduction  and main results}

Let $q$ be a map of $\mathbb{C}$,
\[
 q:x \mapsto x^2+c,
\] 
and
$H_b$ a map of $\mathbb{C}^2$,
\[
 H_{b}:(x,y) \mapsto (x^2+c-by,x),
\]
where parameters $b$, $c$ are complex numbers. They are  respectively the celebrated quadratic map and H\'{e}non map. The latter has Jacobian  $b$, and is a diffeomorphism for nonzero $b$.   
If $b=0$,  $H_{b}$ maps all of $\mathbb{C}^2$ to the graph $x=y^2+c$, hence reduces  to the quadratic map. 

Let $K^+(H_b)$ denote those points whose forward orbits are bounded under iteration of $H_b$,  and $J^+(H_b)$  the boundary of $K^+(H_b)$ in $\mathbb{C}^2$. Define $K^-(H_b)$ and $J^-(H_b)$ analogously by backward iteration of $H_b$.  Define the {\it filled Julia set} by $K(H_b):=K^+(H_b)\cap K^-(H_b)$, and  the {\it Julia set} $J(H_b):=J^+(H_b)\cap J^-(H_b)$.    Note that the Julia set is an invariant set.  (For a map $f$, we call a set $S$   {\it invariant} under $f$ if $f^{-1}(S)=S$.) For the quadratic map, we only consider forward iteration, and in this case  $K(q):=K^+(q)$, $J(q):=J^+(q)$.

Let $\sigma_{q}$ stand for the shift map  in the projective limit  space $\underleftarrow{\lim}(J(q),q)$ induced by $q$. If  $q$ is (uniformly) expanding on its Julia set, i.e. $q$ is hyperbolic, using a perturbation scheme, Hubbard and Oberste-Vorth \cite{HO-V1995} (also see Ishii and Smillie \cite{IS2010}, Radu and Tanase \cite{RT2016}) established the following result. 

\begin{theorem}\label{HO-VFSBS}
If $q$ is hyperbolic and $|b|\not=0$ is sufficiently small,  then the restriction  ${H}_{b}|_{J(H_{b})}$ of $H_{b}$ to its Julia set $J(H_{b})$ is topologically conjugate to  $\sigma_{q}|_{\underleftarrow{\lim}(J(q),q)}$. 
\end{theorem} 

We remark that what they proved is a more general result: $q$ can be generalised to be any  
monic  polynomial  $p:\mathbb{C}\to\mathbb{C}$  of degree  at least two and $H_b$ can be replaced by a generalised  H\'{e}non map 
$f_{p,b}: \mathbb{C}^2\to\mathbb{C}^2, \quad (x,y) \mapsto (p(x)-by,x)$. Note that Bedford and Smillie \cite{BS1992}, Foraeness and Sibony \cite{FS1992},  Hubbard and Oberste-Vorth \cite{HO-V1995},  as well as   Ishii and Smillie \cite{Ishi2008, IS2010} showed that the Julia set $J(f_{p,b})$  is  a (uniformly) hyperbolic set for sufficiently small $|b|$ (thus, the Julia set described in Theorem \ref{HO-VFSBS} is a hyperbolic set).
 
When the Julia set of the H\'{e}non map is a solenoid, Mummert in \cite{Mumm2008} also obtained  a result the same  as Theorem \ref{HO-VFSBS}.
A compact set $S$ is called a {\it solenoid (of degree $2$)} if it is homeomorphic to $\underleftarrow{\lim}(\mathbb{S},x\mapsto x^2)$. If a solenoid $S$ is a hyperbolic invariant set of a diffeomorphism $f$, and $f|_S$ is topologically conjugate to the shift map on $\underleftarrow{\lim}(\mathbb{S}, x\mapsto x^2)$, then we say $f|_S$ 
is a {\it solenoid map}. 
 The subset of parameter space $(c,b)$ for which $H_{b}|_{J(H_{b})}$ is a solenoid map is called the {\it solenoid locus}.
 The following region  
\begin{equation}
 1+2 \sqrt{(1-|b|)^2-4|c|}-4|c|-6|b|-3|b|^2-2|b|\sqrt{(1+|b|)^2+4|c|}>0, \quad c\in \mathcal{M}_0,\quad b\not=0 \label{solenoidMummert}
\end{equation}
that is contained in the solenoid locus was obtained in \cite{Mumm2008}. Here, $\mathcal{M}_0$ is the main cardioid of the Mandelbrot set, i.e. those $c$'s such that $|1-\sqrt{1-4c}|<1$.
 As mentioned  by Mummert  \cite{Mumm2008}, Ishii's work \cite{Ishi2008}  yielded that the  region 
\begin{equation}
 (1-|1- \sqrt{1-4c}|)^2-2|b|-2|b|^2-2|b|\sqrt{(1+|b|)^2+4|c|}>0, \quad c\in \mathcal{M}_0, \quad b\not=0 \label{solenoidIshii}
\end{equation}
is also contained in the solenoid locus. Notice that, for $b=0$, the inequality in \eqref{solenoidIshii} gives $c\in \mathcal{M}_0$, whereas the inequality in \eqref{solenoidMummert} gives $|c|<1/4$. For $c=0$, inequalities (\ref{solenoidIshii}) and \eqref{solenoidMummert} yield $|b|<(\sqrt{2}-1)/2\approx 0.207$ and $|b|<2\sqrt{10}/5-1\approx 0.265$, respectively. Therefore, the parameter region determined by (\ref{solenoidMummert}) does not contain nor is contained by the one determined by (\ref{solenoidIshii}). 

Mummert's method in \cite{Mumm2008} is based on the contraction mapping theorem. Briefly, his method is as follows:   A fixed point of an operator $(x_i)_{i\in\Z}\mapsto (x_i^\prime)_{i\in\Z}$ defined   by 
\begin{equation}
x_i^\prime:=\frac{1}{2}
\left(  x_i+\frac{ x_{i+1}-c+bx_{i-1} }
                      {x_i}
     \right), \qquad i\in\Z, \label{Moperator}
\end{equation}
acting on the complex Banach space $l_\infty (\Z,\C)$ of bounded sequences with the supremum norm
 gives rise to a {bounded orbit} of $H_{b}$.  When $|b|\not= 0$ is small, apply the operator to pseudo-orbits corresponding to points in the projective limit space of $J(q)$ to get genuine orbits of $H_{b}$. (He called this ``shadowing''.) The shadowing process can be continued as long as $\inf_{ (x,y)\in J(H_{b})}
         \min\{ |x|,|y| \}>(1+|b|)/2$. 
 
As pointed out in \cite{Mumm2008}, the operator \eqref{Moperator} is a local linearization of the Sterling-Meiss operator on $l_\infty(\Z, \C)$, $(x_i)_{i\in\Z}\mapsto (x_i^\prime)_{i\in\Z}$ with
\[ x_i^\prime :=s_i \sqrt{x_{i+1}-c+bx_{i-1}}
\]
and $s_i\in\{\pm 1\}$. Sterling and Meiss \cite{SM1998} introduced the above operator to the real H\'{e}non map based on the concept of anti-integrable limit. (We refer the readers to  \cite{AA1990, Chen2010,  LM2006, MM1992, Qin2001} for an account of the anti-integrability.) The fact that this operator is a local contraction in the Euclidean metric enabled them to obtain the Devaney-Nitecki region \cite{DN1979} for the horseshoe locus. Mummert showed that this operator in fact is a local contraction with respect to the Kobayashi metric, thus showed that the Julia set of the complex H\'{e}non map $H_b$ is a horseshoe in the parameter region $|c|>2(1+|b|)^2$ and $b\not= 0$ (also see \cite{Ishi2008, MNTU2000, O-V1987} for obtaining the same region by other approaches).  

The main purpose of this paper is of two folds. One is to reinterpret Mummert's ``shadowing" result in the language of implicit function theorem. The other is to use this shadowing idea to provide an alternative proof of Theorem \ref{HO-VFSBS}. 
Note that the method used in \cite{RT2016} to prove Theorem \ref{HO-VFSBS} is also a contraction mapping theorem applied to a graph-transform operator, while the proofs
 in \cite{HO-V1995, Ishi2008} rely on crossed mappings that are two dimensional analogues of expanding maps.

Our main results are stated in Theorems \ref{mainthm} and \ref{mainthm2}.

Given a compact invariant subset $\Lambda$ of the Julia set $J(q)$, it  induces a compact set, say 
  $\mathcal{A}$, in $\C^{2}$ by
\begin{equation}
 \mathcal{A}:=\bigcup_{x \in \Lambda} (q(x), x).  \label{AfromLambda}
 \end{equation}
Let the mapping $\Lambda\mapsto\mathcal{A}$ be defined by $g$. Clearly, the domain of $g$ can be extended to $\mathbb{C}$, and
\[ H_0\circ g=g\circ q.
\]
Furthermore, $g$ is injective and continuous.

\begin{prop} \label{AH&Lambda} 
~{} \begin{enumerate}
\item $\mathcal{A}$ is an invariant set of the restriction $H_{0}|_{\mathcal{A}}$ of $H_{0}$ to $\mathcal{A}$.
\item  $g$ is a topological conjugacy between $\Lambda$ and $\mathcal{A}$, i.e. the diagram below commutes
 \[
  \begin{matrix}
  \Lambda  & \mapright{q} &  \Lambda  \cr
         \mapldown{g}& &\maprdown{g} \cr
             \mathcal{A} & \mapright{H_{0}}& \mathcal{A}.
  \end{matrix}
\]     
\end{enumerate}
 \end{prop}

We postpone all proofs of theorems in this section until Section \ref{sec:hetero}.

\begin{theorem}\label{mainthm}
 Assume $\Lambda$ is an expanding compact invariant set of $q$. There exists $\epsilon_0>0$ such  that if 
$0<|b|<\epsilon_0$,  then there exists a  hyperbolic compact invariant set $\mathcal{A}_b\subset\mathbb{C}^{2}$ for $H_b$ and homeomorphisms  $\Psi:\underleftarrow{\lim}(\Lambda,q)\to \underleftarrow{\lim}(\mathcal{A},H_{0})$,  $\Upsilon_b:\underleftarrow{\lim}(\mathcal{A},H_{0})\to \mathcal{A}_b$ such that the following diagram commutes 
 \[
  \begin{matrix}
\underleftarrow{\lim}(\Lambda,q)  & \mapright{\sigma_q} &  \underleftarrow{\lim}(\Lambda,q) \cr
         \mapldown{\Psi}& &\maprdown{\Psi} \cr
         \underleftarrow{\lim}(\mathcal{A},H_{0})  & \mapright{\sigma_{\dag}} &  \underleftarrow{\lim}(\mathcal{A},H_{0}) \cr
         \mapldown{\Upsilon_b}& &\maprdown{\Upsilon_b} \cr
             \mathcal{A}_b & \mapright{H_b}& \mathcal{A}_b,  \tag{*} \label{commute_star}
  \end{matrix}
\]
where  $\sigma_{\dag}$ stands for the shift in the projective limit space induced by  $H_{0}$.
\end{theorem}

Recall the definition of holomorphic motion (see  \cite{Jons1998, Lyub1983, MSS1983} for example). Let $b_0\in\C$ be a base point, $\D(b_0, r)\subset \C$ be the open disc of radius $r>0$ centered at $b_0$, and $S$ a subset of $\C^k$, $k\ge 1$. A holomorphic motion of $S$ parameterised by $\D(b_0,r)$ is a continuous map $\Theta:\D(b_0,r) \times S\to \C^k$ such that
\begin{itemize}
 \item $\Theta(b_0, w)=w$ for any $w\in S$.
\item For any $b\in \D(b_0,r)$, the mapping $w\mapsto \Theta(b, w)$ is injective on $S$.
\item For any $w\in S$, the mapping $b\mapsto\Theta(b,w)$ is holomorphic  on $\D(b_0,r)$.
\end{itemize}

\begin{cor} \label{thm:RK278} 
As $b$ varies, the set $\mathcal{A}_b$ in Theorem \ref{mainthm} forms a holomorphic motion. When  $b\to 0$, $\mathcal{A}_b$
converges to $\mathcal{A}$ in the Hausdorff topology,   and   the dynamics degenerate to   $H_0$ on  $\mathcal{A}$. 
 \end{cor}

A converse statement of Theorem \ref{mainthm} is given in Theorem \ref{mainthm2}, which concerns the  opposite direction: 
If $\mathcal{B}_{b_0}$ is a compact hyperbolic invariant set of $H_{b_0}$ for some $b_0\not=0$,  then $\mathcal{B}_{b_0}$ moves holomorphically and 
\begin{equation}
\mathcal{B}_b:=\Theta(b, \mathcal{B}_{b_0}) \label{Aholob}
\end{equation}
is hyperbolic for $H_b$ when $b$ is sufficiently close to $b_0$ (see \cite{Jons1998}).
We have 
\begin{theorem}\label{mainthm2}
Let  $\mathcal{B}_{b_0}\subseteq J(H_{b_0})$  and 
$\mathcal{B}_b$ be described as \eqref{Aholob}. 
There exists $\epsilon_1>0$ such that if  $0<|b_0|<\epsilon_1$, then there exists an expanding compact invariant  set  $\Lambda^\prime \subseteq J(q)$ for $q$ such that $\mathcal{B}_b\to\mathcal{B}=g(\Lambda^\prime)$ in the Hausdorff topology as $b\to 0$.  
 \end{theorem}

A consequence of Theorem \ref{mainthm2} is the following: If $J(H_b)$ is hyperbolic for sufficiently small $|b|$, then it converges to a subset of $g(J(q))$ when $b$ tends to zero. As a result, if $\Lambda=J(q)$ in Theorem \ref{mainthm} and if the set $\mathcal{A}_b$  is a proper subset of $J(H_b)$ for  $b\not=0$, then $\mathcal{A}_b$ would converge to a proper subset of $\mathcal{A}=g(J(q))$ as $b$ approaches zero, contradicting to  Corollary  \ref{thm:RK278}. Hence, we conclude:
\begin{cor}
  Assume the circumstance of Theorem \ref{mainthm}. We have  $\mathcal{A}_b=J(H_b)$ if and only if $\Lambda= J(q)$.
\end{cor}

Combining Theorem \ref{mainthm} with the corollary above, we arrive at a proof of Theorem \ref{HO-VFSBS}. 

We organize this paper as follows.
 In Section \ref{sec:hyperbolicity}, we review some equivalent definitions of hyperbolicity. Our proofs of Theorems \ref{mainthm} and \ref{mainthm2} rely on the implicit function theorem, thus in Section \ref{sec:invertible} we define  linear operators in the spaces of bounded sequences and study their invertibility. 
Section \ref{sec:hetero} is devoted to proving Theorems \ref{mainthm} and \ref{mainthm2}. 
A crucial part of our approach will be based on a theorem that initially is developed for the theory of anti-integrable limit.    For the sake of the completeness of our exposition, we include in the Appendix the theorem and a slighted modified proof by adding the uniform discreteness for it  .

\section{Uniform hyperbolicity} \label{sec:hyperbolicity}

In this section, we review some  definitions of the  hyperbolicity. 
 Recall first the definition of cone fields. 
 Given $k\ge 2$, $c>0$, a norm $\|\cdot\|$ in $\C^k$, and a splitting of the tangent space $T_z\mathbb{C}^k=V_z^s\oplus V_z^u$ at each point $z$, the {\it stable} and {\it unstable} $c$-{\it cones} are defined by 
 \begin{eqnarray*}
       \mathcal{C}_z^s=\{w+v \in V_z^s\oplus V_z^u| ~ \|v\|\le c\|w\|\}\supset V_z^s, \\
      \mathcal{C}_z^u=\{w+v \in V_z^s \oplus V_z^u| ~ \|w\|\le c\|v\|\}\supset V_z^u,
      \end{eqnarray*}
respectively.      A compact invariant set $S$ of a diffeomorphism $f$ of $\C^k$ is said to satisfy the {\it cone criterion} if  there exist $0<\lambda_s<1<\lambda_u$ such that for every $z\in S$ there is a splitting of the tangent space and families of cones   
       associated with the splitting such that 
\[ Df(z)\mathcal{C}_z^u\subset\interior\mathcal{C}_{f(z)}^u, \qquad Df(z)^{-1}\mathcal{C}_{f(z)}^s\subset\interior\mathcal{C}_{z}^s,
\]
$\|Df(z)\zeta \|\ge \lambda_u \|\zeta\|$ for $\zeta\in\mathcal{C}_z^u$, and that  $\|Df(z)^{-1}\zeta\|\ge \lambda_s^{-1} \|\zeta\|$ for $\zeta\in\mathcal{C}_{f(z)}^s$ and some norm $\|\cdot\|$ in $\C^k$.

  A sequence of invertible linear maps $A_i:\mathbb{C}^k\to\mathbb{C}^k$, $i\in\mathbb{Z}$ and $k\ge 2$, is said to admit a {\it $(\lambda_s,\lambda_u)$-splitting} if there exist direct sum decompositions $\mathbb{C}^k=E_i^s\oplus E_i^u$ such that $A_iE_i^{s/u}=E_{i+1}^{s/u}$ and there are constants $C_s$, $C_u>0$,  $0<\lambda_s <1< \lambda_u$ such that  $\|A_{i+n-1}\cdots A_{i+1}A_i|_{E^s_i}\|\le C_s \lambda_s^n$ and  $\|A_{i-n+1}^{-1}\cdots A_{i-1}^{-1}A_i^{-1}|_{E_{i+1}^u}\|\le C_u \lambda_u^{-n}$ for some norm $\|\cdot\|$ in $\C^k$, all $i\in\mathbb{Z}$ and  $n\ge 1$.    
 Such a sequence of maps $A_i$ is said to have an {\it exponential dichotomy} on $\mathbb{Z}$ in \cite{Palm2000}. The subspaces $E_i^s$ and $E_i^u$ are called the {\it stable} and {\it unstable subspaces}. We refer to the complex dimension of the unstable subspace as the {\it unstable dimension}.  
 
\begin{defn} \label{defnhyp} \rm
  Let $f:\mathbb{C}^k\to\mathbb{C}^k$, $k\ge 2$, be a diffeomorphism. Suppose that $S$ is a compact $f$-invariant set. Then $S$ is called a {\it (uniformly) hyperbolic} set for $f$ if for any $z\in S$ the sequence of differentials $Df(f^i(z)): T_{f^i(z)}\mathbb{C}^k\to T_{f^{i+1}(z)}\mathbb{C}^k$, $i \in\mathbb{Z}$, admits a $(\lambda_s,\lambda_u)$-splitting and the splitting varies continuously with the point $z$.
\end{defn} 

The notion of hyperbolicity for a map of $\C$ is associated with the expansion on its Julia set. 
\begin{defn} \label{defnexp} \rm
An invariant set $S$ for a differentiable map $f$ of $\mathbb{C}$ is {\it (uniformly)  expanding} if there exist $C_u>0$ and $\lambda_u>1$ such that for all $z\in S$ and $n\ge 1$ we have $\|Df^n(z)\|\ge C_u\lambda_u^n$ for some norm $\|\cdot\|$ in $\C$.
  \end{defn}

Given a  differentiable map  $f$ of $\mathbb{C}^k$, $k\ge 2$, we introduce a map $F(\cdot; f)$ on the Banach space $l_\infty (\mathbb{Z},\mathbb{C}^k):=\{  {\bf z}|\ {\bf z}=(\ldots,z_{-1},z_0,z_1,\ldots),~z_i\in\mathbb{C}^k,\ \mbox{bounded}\}$ of bounded sequences with the supremum 
norm $\|\cdot\|_\infty$:
\begin{eqnarray}
  F(\cdot; f):l_\infty (\mathbb{Z},\mathbb{C}^k) &\to& l_\infty (\mathbb{Z},\mathbb{C}^k), \nonumber \\
  {\bf z} &\mapsto& F({\bf z}; f)=\left(F({\bf z}; f)_{i}\right)_{i\in\mathbb{Z}} \label{MapFfBanach}
\end{eqnarray}
with $F({\bf z}; f)_i= z_{i+1}-f(z_i)$. It is readily to see that ${\bf z}$ is a bounded orbit of $f$ if and only if it solves $F({\bf z}; f)=0$.
  Certainly, $F(\cdot; f)$ is a differentiable map, with its derivative at ${\bf z}$ a bounded linear operator:
\begin{eqnarray*}
  DF({\bf z}; f):l_\infty (\mathbb{Z},\mathbb{C}^k) &\to&  l_\infty (\mathbb{Z},\mathbb{C}^k), \\
   (\xi_i)_{i\in\Z}=\boldsymbol{\xi} &\mapsto&  DF({\bf z}; f)\boldsymbol{\xi}=\left(\sum_{j\in\mathbb{Z}} D_{z_j}F({\bf z}; f)_i~\xi_j\right)_{i\in\mathbb{Z}}=(\eta_i)_{i\in\Z}=\boldsymbol{\eta}.
\end{eqnarray*}
 
In addition to Definitions \ref{defnhyp} and \ref{defnexp}, the notion of hyperbolicity in several different situations   can be characterized by different ways, for instance, by the phonon gap \cite{AMB92}, by the invertibility of an operator in the tangent bundle of an orbit  \cite{Lanf1985}, by the disjointness of the Mather spectrum from the unit circle \cite{HK2002, Math1968}, by the cone criterion (e.g. \cite{NP1973}), and by the exponential dichotomy (e.g. \cite{Palm2000}). Among them, the following four  statements are particularly useful in our approach.  

\begin{theorem} \label{hypequiv}
 Let $S$ be a compact invariant set for a diffeomorphism  $f$ of $\mathbb{C}^k$, $k\ge 2$, and ${\bf z}=(z_i)_{i\in\Z}$ an orbit of $f$ with $z_i=f^i(z_0)$ and $z_0\in S$.  The following statements are equivalent.
\begin{enumerate}
\item $S$ is  hyperbolic. 
\item $S$ satisfies the cone criterion.
\item For every orbit  ${\bf z}$, the derivative $DF({\bf z}; f)$ is an invertible linear operator of $l_\infty(\mathbb{Z},\mathbb{C}^k)$ with $\|DF({\bf z}; f)^{-1}\|_\infty$ bounded uniformly in $\mathbf{z}$. 
\item For any given $\boldsymbol{\eta}=(\eta_i)_{i\in\Z}\in l_\infty(\mathbb{Z},\mathbb{C}^k)$ with $\|\boldsymbol{\eta}\|_{\infty}=1$ and any orbit ${\bf z}$, the recurrence relation 
  \begin{equation}
   \xi_{i+1}-Df(z_i)\xi_i=\eta_i,\quad i\in\mathbb{Z}, \label{hypequivnonhomo}
  \end{equation}
    has a unique solution $\boldsymbol{\xi}=(\xi_i)_{i\in\mathbb{Z}}\in l_\infty(\mathbb{Z},\mathbb{C}^k)$ such that $\|\boldsymbol{\xi}\|_\infty$ is uniformly bounded in ${\bf z}$. 
\end{enumerate}
\end{theorem}

\begin{rk} \rm \label{InverseNorm}
In the theorem above, if the tangent bundle $T_S\mathbb{C}^k$ of the hyperbolic set $S$  splits into a direct sum of two invariant sub-bundles for the tangent map $Df$ with constants $C_s$, $C_u$, and a $(\lambda_s,\lambda_u)$-splitting, then  we have  $\|DF({\bf z}; f)^{-1}\|_\infty \le C_s(1-\lambda_s)^{-1}+C_u(\lambda_u-1)^{-1}$ for every ${\bf z}=(f^i(z_0))_{i\in\Z}$ and $z_0\in S$ (see e.g. \cite{Palm2000}).
\end{rk}   

For a non-invertible holomorphic  map $f$ of $\mathbb{C}$, we can rewrite it as a map $F(\cdot; f)$ on the Banach space $l_\infty (\mathbb{Z}^+,\mathbb{C}):=\{  {\bf x}^+|\ {\bf x}^+ =(x_0,x_1,\ldots), ~z_i\in\mathbb{C},\ \mbox{bounded}\}$ of bounded sequences with the supremum norm $\|\cdot\|_\infty$:
\begin{eqnarray}
  F(\cdot; f):l_\infty (\mathbb{Z}^+,\mathbb{C}) &\to& l_\infty (\mathbb{Z}^+,\mathbb{C}), \nonumber\\
  {\bf x}^+ &\mapsto& F({\bf x}^+; f)=\left(F({\bf x}^+; f)_i\right)_{i\ge 0} \label{MapFqBanach}
\end{eqnarray}
with $F({\bf x}^+; f)_i= x_{i+1}-f(x_i)$ for all $i\ge 0$. Then ${\bf x}^+$ is a bounded (forward) orbit of $f$ if and only if it solves $F({\bf x}^+; f)=0$.
  The map $F(\cdot; f)$ is  differentiable with its derivative at ${\bf x}^+$ a bounded linear operator:
\begin{eqnarray*}
  DF({\bf x}^+; f):l_\infty (\mathbb{Z}^+,\mathbb{C}) &\to&  l_\infty (\mathbb{Z}^+,\mathbb{C}), \\
  (\xi_i)_{i\ge 0}= \boldsymbol{\xi}^+ &\mapsto&  DF({\bf x}^+; f)\boldsymbol{\xi}^+=\left(\sum_{j=0}^\infty D_{z_j}F({\bf x^+}; f)_i~\xi_j\right)_{i\ge 0}=(\eta_i)_{i\ge 0}=\boldsymbol{\eta}^+.
\end{eqnarray*}
Theorem \ref{hypequiv} has a counterpart for non-invertible maps of $\C$, as follows.
\begin{theorem} \label{thm1}
 Let $S$ be a compact  invariant set for a holomorphic map $f:\mathbb{C}\to\mathbb{C}$, then the following  statements are equivalent. 
\begin{enumerate}
\item $S$ is expanding. 
\item For any orbit ${\bf x}^+=(x_i)_{i\ge 0}$ of $f$ with $x_i=f^i(x_0)$ and $x_0\in S$, the derivative $DF({\bf x}^+; f)$ is an invertible linear operator of $l_\infty(\mathbb{Z}^+,\mathbb{C})$ with $\|DF({\bf x}^+; f)^{-1}\|_\infty$ bounded uniformly in $\mathbf{x}^+$. 
\item  For any given $\boldsymbol{\eta}^+=(\eta_i)_{i\ge 0}\in l_\infty(\mathbb{Z}^+,\mathbb{C})$ and any orbit $(x_i)_{i\ge 0}$ with $x_0\in S$, the recurrence relation
    \begin{equation}
      \xi_{i+1}-Df(x_i)\xi_i=\eta_i,\quad i\ge 0, \label{nonauto}
    \end{equation}
    has a unique bounded solution $\boldsymbol{\xi}^+=(\xi_i)_{i\ge 0}\in l_\infty(\mathbb{Z}^+,\mathbb{C})$
\end{enumerate}
\end{theorem} 

\begin{rk} \rm \label{InverseNormL=0}
In the theorem above, if the tangent map has an expansion rate $\lambda_u$ with a constant $C_u$ on $T_S\mathbb{C}$, then  the proof given in \cite{Palm2000} for Remark \ref{InverseNorm} also works to give  $\|DF({\bf x}; f)^{-1}\|_\infty \le C_u(\lambda_u-1)^{-1}$ for every ${\bf x}=(f^i(x_0))_{i\ge 0}$ and $x_0\in S$ 
\end{rk}   

 \begin{rk} \label{rk:necessary} \rm 
The statement (iv) of Theorem \ref{hypequiv}  implies that the only bounded homogeneous solution of (\ref{hypequivnonhomo})  is necessarily trivial.
\end{rk}
The necessary condition in the above remark, however, turns out to be sufficient when $S$ is chain recurrent. The following result is due to Churchill {\it et al}. \cite{CFS1977}, Sacker and Sell \cite{SS1974}. 
\begin{theorem} \label{thm:CFS}
Assume  $S$ is a compact chain recurrent invariant set of  a diffeomorphism $f$ of $\C^k$, $k\ge 2$. Then $S$ is  hyperbolic  if and only if the recurrence relation
 \begin{equation}
   \xi_{i+1}-Df(z_i)~ \xi_i=0,\quad i\in\mathbb{Z}, \label{hypequivhomo}
 \end{equation}
 has no non-trivial bounded solutions  $(\xi_i)_{i\in\Z}$ for any orbit $(z_i)_{i\in\Z}$ with $z_0\in S$. 
\end{theorem}
A compact $f$-invariant set $S$ is called {\it quasi-hyperbolic} in \cite{CFS1977} if  the recurrence relation \eqref{hypequivhomo}  has no non-trivial bounded solutions. Based on Theorem \ref{thm:CFS}, Arai \cite{Arai2007} developed a rigorous computational method to prove the quasi-hyperbolicity, and obtained the hyperbolic plateaus of the real H\'{e}non maps.  

The following is the  version of quasi-hyperbolicity for non-invertible maps of $\C$ (see  \cite{Chen2008}). 
\begin{theorem} 
Assume  $S$ is a compact invariant set of  a holomorphic map $f$ of $\C$. Then $S$ is expanding  if and only if the recurrence relation
 \begin{equation}
   \xi_{i+1}-Df(x_i)~ \xi_i=0,\quad i\ge 0 , \label{hypequivhomoC}
 \end{equation}
 has no non-trivial bounded solutions  $(\xi_i)_{i\ge 0}$ for any orbit $(x_i)_{i\ge 0}$ with $x_0\in S$.
\end{theorem}

We end this section by a remark that the terminology `quasi-hyperbolicity' means differently in \cite{BS2002}. There, a diffeomorphism of $\C^2$ is called {\it quasi-expanding} if a certain canonical metric is  expanding on unstable subspaces of saddle periodic points, and is called {\it quasi-contracting} if its inverse is quasi-expanding. The diffeomorphism is quasi-hyperbolic if it is both quasi-expanding and quasi-contracting. They showed that the quasi-expansion may be viewed as a two-dimensional analogue of semi-hyperbolicity for one-dimensional polynomial maps.

\section{Invertibility of linear operators} \label{sec:invertible}

Recall   the {\it projective limit} (or called the {\it inverse limit}) {\it space}  for a  continuous map $f$ of $\C^k$, $k\ge 1$. It is defined by the following subset of $(\C^k)^{\mathbb{Z}^-}$:
\[   \underleftarrow{\lim} (\C^k,f):=\{(\ldots,z_{-1},z_0)\in (\C^k)^{\mathbb{Z}^-}|~ f(z_{i-1})=z_i  ~ \forall ~ i\le 0\}.
\]
A point in this space is a point $z_0$ and a history of $z_0$ under iteration of $f$. 
The map $f$ induces a shift map
\begin{eqnarray}
  \sigma_{f}:\underleftarrow{\lim} (\C^k,f) \to \underleftarrow{\lim} (\C^k,f), \quad
  (\ldots,z_{-2},z_{-1},z_0) \mapsto (\ldots,z_{-1},z_0,f(z_0)). \label{sigmaf}
\end{eqnarray}
Obviously, $\sigma_{f}$ is a bijection and 
\[
  \sigma_{f}^{-1}(\ldots,z_{-2},z_{-1},z_0)=(\ldots,z_{-3},z_{-2},z_{-1}).
\]
Endow the projective limit space with the product topology. 

 Subsequently, $\sigma_f$ is a homeomorphism.
Let $S$ be a bounded $f$-invariant subset of $\C^k$, we use
\[   \underleftarrow{\lim} (S,f):=\{(\ldots,z_{-1},z_0)
 \in\underleftarrow{\lim} (\C^k,f)|~ z_i\in S ~ \forall ~ i\le 0\}.
\]

Define 
\begin{eqnarray*}
 \Sigma(S,f)&:=&\{(\ldots,z_{-1},z_0,z_1,\ldots)\in S^{\mathbb{Z}}|~ f(z_{i})=z_{i+1}~ \forall i\in\mathbb{Z}\},  \\
 \Sigma^+(S,f)&:=& \{(z_0, z_1, z_2, \ldots)\in S^{\mathbb{Z}^+}|~ f(z_{i})=z_{i+1}~\forall i\ge 0\}.
 \end{eqnarray*}
 Notice that $\Sigma(S,f)$ with the product topology is homeomorphic to $\underleftarrow{\lim}(S,f)$.
For any ${\bf z}=(\ldots, z_{-1}, z_0, z_1, \ldots)\in (\C^k)^\Z$, define the truncation $\pi^+({\bf z})$ of ${\bf z}$ to be an element of $(\C^k)^{\mathbb{Z}^+}$ by 
\[ \pi^+({\bf z})=(z_0,z_1,\ldots).
\] 
 
Concerning  bounded orbits of the maps $q$, $H_0$ and $H_b$, with the notation and the function $F$ defined in \eqref{MapFfBanach} or \eqref{MapFqBanach} in previous section, explicitly we shall consider respectively the following three maps:
\[ F({\bf x}^+; q)_i=x_{i+1}-q(x_i), \quad  i\ge 0, \quad
\mbox{with} ~ {\bf x}^+=(x_i)_{i\ge 0}\in l_\infty(\Z^+, \C), 
\]
\begin{equation}
 F({\bf w}^\dag; H_0)_i=z_{i+1}-H_0(z_i), \quad  i\in\Z, \quad\mbox{with}~ {\bf w}^\dag=(z_i)_{i\in\Z}\in l_\infty(\Z, \C^{2}),  \label{11-1}
\end{equation}
and 
\begin{equation}
 F({\bf w}; H_b)_i=z_{i+1}-H_b(z_i), \quad  i\in\Z, \quad \mbox{with}~{\bf w}=(z_i)_{i\in\Z}\in l_\infty(\Z, \C^{2}).\label{defnFwHi}
\end{equation}

The Lemma  below is the key lemma of this paper. 

\begin{la} \label{keylemma}
 Suppose that $\mathcal{A}\subset\C^{2}$ is defined as (\ref{AfromLambda}).
 The linear operator $DF({\bf w}^\dag;H_{0}):l_\infty(\Z,\C^{2})\to l_\infty(\Z,\C^{2})$ is invertible for every ${\bf w}^\dag\in\Sigma(\mathcal{A},H_{0})$ if and only if
$\Lambda$ is  expanding for $q$.
\end{la}
  \proof
Suppose  
${\bf x}=(x_i)_{i\in\Z}\in\Sigma(\Lambda,q)$. Let ${\bf w}^\dag=(w_i)_{i\in\Z}=(x_{1+i},  x_{i})_{i\in\Z}$. 
 We have  ${\bf w}^\dag\in \Sigma(\mathcal{A},H_{0})$.

  Theorem \ref{thm1} tells that $\Lambda$ is  expanding for $q$ if and only if the recurrence relation 
  \[
   \xi_{i+1}-Dq(x_i)\xi_{i}= {\eta}_{i}, \quad\forall i\ge 0,
   \]
  has a unique solution  $(\xi_{i})_{i\ge 0}\in l_\infty(\Z^+,\C)$ for any given  $({\eta}_{i})_{i\ge 0}\in l_\infty(\Z^+,\C)$. 
The solution depends on the orbit ${\bf x}^+=(x)_{i \ge 0}$ of $q$  
and on the sequence ${\boldsymbol{\eta}}^+=({\eta}_{i})_{i \ge 0}$, thus we write it as 
\[ (\xi_{i})_{i\ge 0}=:\boldsymbol{\xi}^+=\boldsymbol{\xi}^+({\boldsymbol{\eta}}^+; ~{\bf x}^+)=(\boldsymbol{\xi}^+({\boldsymbol{\eta}}^+; ~{\bf x}^+)_i)_{i\ge 0}.
\]
Consequently, 
\[
   \boldsymbol{\xi} = \boldsymbol{\xi}({\boldsymbol{\eta}};~{\bf x})=(\ldots,  \xi_{-1}, \xi_{0}, \xi_{1}, \ldots) \quad \in l_\infty(\Z,\C) 
 \]
obtained by  
\begin{equation}
 \xi_{i}=\begin{cases}
                     \boldsymbol{\xi}^+\left(\pi^+\circ\sigma^i ({\boldsymbol\eta}); ~\pi^+\circ\sigma^i ({\bf x})\right)_0 & \mbox{for} ~ i\le 0 \\
                     \boldsymbol{\xi}^+\left(\pi^+ ({\boldsymbol\eta}); ~\pi^+ ({\bf x})\right)_i & \mbox{for}~ i\ge 0 
                   \end{cases}   \label{xiis!}
\end{equation}
   is a bounded solution of the recurrence relation
  \begin{equation}
   \xi_{i+1}-Dq(x_i)\xi_{i}= {\eta}_{i}, \quad\forall i\in\Z,  \label{eq30+1}
   \end{equation}
for any given ${\boldsymbol{\eta}}=({\eta}_{i})_{i\in\Z}\in l_\infty(\Z,\C)$, where ${\bf x}=(x_{i})_{i \in\Z}$ satisfies $x_i\in q^{-1}(x_{i+1})$.
 It is straightforward to see that such a solution $\boldsymbol{\xi}$  given by (\ref{xiis!}) is unique. 

To show that $DF({\bf w}^\dag;H_{0})$ is invertible, it amounts to showing that,  for any given bounded sequence   $\tilde{\boldsymbol{\eta}}=({\eta}_{i},\tilde\eta_{i})_{i\in\Z}\in l_\infty(\Z, \C^{2})$,  the recurrence 
  relations
 \begin{eqnarray}
  \xi_{i+1}-Dq( x_{i})\xi_{i} &=& \eta_{i}, \label{la31m_1} \\
  \tilde{\xi}_{i+1}-\xi_{i} &=& \tilde{\eta}_{i}, \label{xi(01)eta(1)}
 \end{eqnarray}
 for all $i\in\Z$
  have a unique bounded solution $\tilde{\boldsymbol{\xi}}=({\xi}_{i},\tilde\xi_{i})_{i\in\Z}\in l_\infty(\Z, \C^2)$. 
 Now, equations \eqref{la31m_1} and \eqref{eq30+1} are the same. Once $\xi_{i}$ for every integer $i$ is obtained, \eqref{xi(01)eta(1)} can be solved by
\begin{equation}
 \tilde\xi_{i} = \tilde\eta_{i-1}+\xi_{i-1} \label{xii-1}   
\end{equation}
 for all $i\in\Z$.   This proves the invertibility of $DF({\bf w}^\dag;H_{0})$.

The proof of the other direction of the lemma is easy to see.
  \qed \\

The following result provides an estimate on the size of $DF({\bf w}^\dag; H_0)^{-1}$  (cf. Remarks \ref{InverseNorm} and \ref{InverseNormL=0}).
\begin{prop} \label{prop:InverseBoundHphi}
For every ${\bf w}^\dag\in\Sigma(\mathcal{A},H_0)$,  we have
\[ \|  
         DF(     {\bf w}^\dag; H_0
              )^{-1}
    \|_\infty \le 
              \frac{2C_u}{\lambda_u-1}+1
\]
 if $Dq$ is expanding on $T\Lambda$ with rate $\lambda_u$ and constant $C_u$. 
\end{prop}
\proof
    Let ${\bf x}=(x_i)_{i\in\Z}\in\Sigma(\Lambda,q)$, and  ${\bf w}^\dag=(w_i)_{i\in\Z}=(x_{1+i},x_{i})_{i\in\Z}\in \Sigma(\mathcal{A}, H_0)$. By definition, 
\[ \|DF({\bf w}; H_0)^{-1}\|_\infty=\sup_{\tilde{\boldsymbol{\eta}}}\|\tilde{\boldsymbol{\xi}}(\tilde{\boldsymbol{\eta}})\|_\infty =\sup_{\tilde{\boldsymbol{\eta}}}\sup_{i\in\Z}|(\xi_{i}(\tilde{\boldsymbol{\eta}}), \tilde\xi_{i}(\tilde{\boldsymbol{\eta}}))|
\]
with $\tilde{\boldsymbol{\xi}}=\tilde{\boldsymbol{\xi}}(\tilde{\boldsymbol{\eta}})$ satisfying the recurrence relations \eqref{la31m_1} and \eqref{xi(01)eta(1)}  for $\tilde{\boldsymbol{\eta}}\in l_\infty(\Z,\C^{2})$ and $\|\tilde{\boldsymbol{\eta}}\|_\infty=1$. By the triangular inequality, 
\begin{eqnarray*}
 |(\xi_{i}(\tilde{\boldsymbol{\eta}}), \tilde\xi_{i}(\tilde{\boldsymbol{\eta}}))|  &\le&
          |(\xi_{i}(\tilde{\boldsymbol{\eta}}), 0)|+|(0, \tilde\xi_{i}(\tilde{\boldsymbol{\eta}}))| \\
&\le & |\xi_{i}(\tilde{\boldsymbol{\eta}})|+|\tilde\xi_{i}(\tilde{\boldsymbol{\eta}})| \\  
&\le & |\xi_{i}(\tilde{\boldsymbol\eta})|+|\xi_{i-1}(\tilde{\boldsymbol\eta})|+1 \qquad\mbox{(by equality \eqref{xii-1}).} 
\end{eqnarray*}
 Let ${\bf x}^+ = \pi^+({\bf x})$, and $\boldsymbol{\xi}^+(\pi^+((\eta_i)_{i\in\Z}); \pi^+({\bf x}))=DF({\bf x}^+;q)^{-1}\pi^+((\eta_i)_{i\in\Z})$. Now,
\[
 \xi_{i}((\eta_i)_{i\in\Z})=\begin{cases}
                     \boldsymbol{\xi}^+\left(\pi^+\circ\sigma^i ((\eta_i)_{i\in\Z}); ~\pi^+\circ\sigma^i ({\bf x})\right)_0 & \mbox{for} ~ i\le 0 \\
                     \boldsymbol{\xi}^+\left(\pi^+ ((\eta_i)_{i\in\Z}); ~\pi^+ ({\bf x})\right)_i & \mbox{for}~ i\ge 0 
                   \end{cases}   
\]
as described in \eqref{xiis!}.  Therefore, we have for $\|(\eta_i)_{i\in\Z}\|_\infty=1$ that
\begin{eqnarray*}
 \sup_{\tilde{\boldsymbol\eta}}\sup_{i\in\Z} \left(
                                                                            |\xi_{i}(\tilde{\boldsymbol\eta})|
                                                                                      +|\xi_{i-1}(\tilde{\boldsymbol\eta})|                                                                                                                                \right)
&\le& 2 \sup_{{\bf x}^+\in\Sigma(\Lambda, q)}\sup_{(\eta_i)_{i\in\Z}}\|DF({\bf x}^+;q)^{-1}\pi^+((\eta_i)_{i\in\Z})\|_\infty \\
&= &2 \sup_{{\bf x}^+\in\Sigma(\Lambda, q)}\|DF({\bf x}^+,q)^{-1})\|_\infty.
\end{eqnarray*}
In view of Remark \ref{InverseNormL=0}, this leads to the proposition.
\qed

\section{Topological conjugacy} \label{sec:hetero}

This section is mainly dedicated to proving Theorems \ref{mainthm} and \ref{mainthm2}. Before that, we prove Proposition  \ref{AH&Lambda} first.

\subsection{Proof of Proposition \ref{AH&Lambda}}

(i) Suppose we are given a point $(x, y)\in\mathcal{A}$, then from the definition of $\mathcal{A}$, the point
$y$ belongs to $\Lambda$ and $  q ( y)=x$.
By the $q$-invariance of $\Lambda$, we have $x \in\Lambda$. Now, $ (x,y) \in H_0|_{\mathcal{A}}^{-1} ( q (x), x)$ and $(q(x), x)\in\mathcal{A}$. Hence,  $\mathcal{A}\subset H_0|_{\mathcal A}^{-1}(\mathcal{A})$. 
To show $\mathcal{A}\supset H_0|_{\mathcal{A}}^{-1}(\mathcal{A})$, we have to find a point in $\mathcal{A}$ whose image
under $H_0$ is $(x,y)$. A candidate is $(y, y')$ where $y'$ belongs to $\Lambda$ and satisfies $q (y')=y$. Such $y'$ exists because $y$ belongs to $\Lambda$ and the set $\Lambda$ is invariant under $q$.  
(ii) $g$ is a continuous bijection on the compact set $\Lambda$, hence a homeomorphism.  \qed

\subsection{Proof of Theorem \ref{mainthm}}

The existence of the homeomorphism $\Psi$ with which $\sigma_q|_{\underleftarrow{\lim}(\Lambda, q)}$   is topologically conjugate to $\sigma_\dag|_{\underleftarrow{\lim}(\mathcal{A}, H_0)}$
follows directly from Proposition \ref{AH&Lambda}(ii). The rest of the proof is devoted to the existence of $\mathcal{A}_b$ and of the topological conjugacy $\Upsilon_b$ between $\underleftarrow{\lim}(\mathcal{A}, H_0)$ and $\mathcal{A}_b$. Because the two sets $\underleftarrow{\lim}(\mathcal{A}, H_0)$ and  $\Sigma(\mathcal{A}, H_0)$ are homeomorphic to each other, it is enough to work on the latter set. 

\subsubsection*{Continuation:}

 If $\Lambda$ is a compact  hyperbolic invariant set for $q$, then the linear operator $DF({\bf w}^\dag;H_0)$ is invertible for every ${\bf w}^\dag\in \Sigma(\mathcal{A}, H_0)$  by  Lemma \ref{keylemma}.  Since $F({\bf w}^\dag;H_0)=0$, by virtue of  the implicit function theorem (e.g. \cite{Whit1965}),  there exists $\hat\epsilon=\hat\epsilon({\bf w}^\dag)$ and a unique holomorphic function 
  \begin{equation}
 \theta(\cdot;{\bf w}^\dag):\C\to l_\infty(\Z,\C^{2}), \qquad b\mapsto   \theta(b;{\bf w}^\dag)=(\theta(b;{\bf w}^\dag)_i)_{i\in\Z} \label{thetabw}
  \end{equation}
  such that $F(\theta(b;{\bf w}^\dag);H_b)=0$ and $\theta(0;{\bf w}^\dag)={\bf w}^\dag$ provided $|b|<\hat\epsilon$. The construction of $F$ tells that $\theta(b;{\bf w}^\dag)$ is a bounded orbit of the H\'{e}non map $H_b$.

\subsubsection*{Bijectivity:}

Note that $\Lambda$ is expansive. That is, there exists $\tau>0$ such that for distinct $x$, $\tilde{x}\in\Lambda$, there is some $n\ge 0$ such that $|q^n(x)-q^n(\tilde{x})|>\tau$. This implies that $\Sigma(\Lambda, q)$ viewed as a subset of $l_\infty(\Z,\C)$ is uniformly discrete, so is $\Sigma(\mathcal{A}, H_0)$ as a subset of $l_\infty(\Z,\C^2)$.  (A subset $\Sigma$ of $l_\infty(\Z, \C^{k})$, $k\ge 1$, being {\it uniformly discrete} means that there exists $\tau>0$ such that whenever ${\bf u}$ and ${\bf v}$ are distinct points of $\Sigma$, then $\|{\bf u}-{\bf v}\|_\infty >\tau$.)  By Proposition \ref{prop:InverseBoundHphi}, $\|DF({\bf w}^\dag;H_0)^{-1}\|_\infty$ is bounded above on $\Sigma(\mathcal{A},H_0)$ when $q$ (thus $H_0$) is given. From the definitions of $F({\bf w}^\dag; H_0)$ and $F({\bf w}; H_b)$ (see \eqref{11-1} and \eqref{defnFwHi}), it is clear that 
for any $\gamma>0$ there exist $\lambda_0>0$ and $\delta_0>0$ such that for all ${\bf w}^\dag\in\Sigma(\mathcal{A},H_0)$ we have $\|F({\bf w};H_b)\|_\infty<\gamma$ and $\|DF({\bf w};H_b)-DF({\bf w}^\dag;H_0)\|_\infty<\gamma$ whenever $\|{\bf w}-{\bf w}^\dag\|_\infty< \lambda_0$ and $|b|<\delta_0$. Therefore, by making use of Theorem \ref{la23}(i) in the Appendix, we have $\epsilon_0:=\inf_{{\bf w}^\dag\in\Sigma(\mathcal{A}, H_0)}\hat{\epsilon}({\bf w}^\dag)> 0$,
and the mapping
\[        \Sigma(\mathcal{A}, H_0) \to  \bigcup_{{\bf w}^\dag\in\Sigma(\mathcal{A},H_0)}\theta(b; {\bf w}^\dag)
   \]
is a bijection provided $|b|<\epsilon_0$.

\subsubsection*{Conjugacy:}

  Let the projection ${\bf w}=(\cdots,w_{-1},w_0,w_1,\cdots)\mapsto w_0 \in\mathbb{C}^{2}$ be  denoted by $\pi_0$.  
 We assert that   the  composition of mappings  
  \begin{equation}
{\bf w}^\dag\ \stackrel{\Phi_b}{\longmapsto}\ \theta(b;{\bf w}^\dag)\ \stackrel{\pi_0}{\longmapsto}~ \theta(b;{\bf w}^\dag)_0 \label{Phib}
  \end{equation}
  is a homeomorphism with the product topology, and that  the  following diagram 
\begin{equation}
  \begin{matrix}
        \Sigma (\mathcal{A}, H_0)& \mapright{\sigma_\dag} &   \Sigma (\mathcal{A}, H_0) \cr
         \mapldown{\pi_0\circ\Phi_b}& &\maprdown{\pi_0\circ\Phi_b} \cr
             \mathcal{A}_b& \mapright{H_b}& \mathcal{A}_b
  \end{matrix} \label{diadiacom}
\end{equation}
commutes when $0<|b|<\epsilon_0$, where the set $\mathcal{A}_b$ is invariant under $H_b$ and is defined by 
  \[ \mathcal{A}_b:=\bigcup_{
                                                {\bf w}^\dag\in   \Sigma (\mathcal{A};H_0)
                                            }                                       \pi_0\circ\theta(b;{\bf w}^\dag).
\]
In other words, $\pi_0\circ\Phi_b$ acts as the topological conjugacy.

  The projection  $\pi_0$ certainly is continuous, and it is bijective since the sequence $\theta(b;{\bf w}^\dag)$ is an orbit of $H_b$ and is uniquely determined by the initial points $\theta(b;{\bf w}^\dag)_0$ when $H_b$ is a diffeomorphism. We showed that  $\Phi_b$ is bijective. By Proposition \ref{AH&Lambda}(ii), $\Lambda$ and $\mathcal{A}$ are topologically the same, thus the compactness of $\Lambda$ implies the compactness of  $\Sigma(\mathcal{A}, H_0)$ in the product topology. Hence, we can set $\epsilon_0$ to be smaller if necessary so that the requirement of Theorem \ref{la23}(ii) is fulfilled and 
accordingly $\Phi_b$ is continuous. Then $\pi_0\circ\Phi_b$ is a continuous bijection from a compact space to a Hausdorff space, thus a homeomorphism. 
Note that $q$-invariance of $\Lambda$ implies $\sigma_{\dag}$-invariance of  
$\Sigma(\mathcal{A}, H_0)$, and that 
\begin{eqnarray*}
F({\bf w};H)_{i+1} &=& w_{i+2}-H(w_{i+1}) \\
            &=&\sigma ({\bf w})_{i+1}-H(\sigma ({\bf w})_{i}) \\
          &=& F(\sigma ({\bf w}); H)_i,
\end{eqnarray*}
 thence the commutativity of the diagram \eqref{diadiacom} follows immediately from Theorem \ref{la23}(iii) and from the bijectivity of $\pi_0$. (If there is no confusion, we use $\sigma_\dag$ to represent the shift operator in both $\Sigma(\mathcal{A}, H_0)$ and $\underleftarrow{\lim}(\mathcal{A}, H_0)$.)

\subsubsection*{Hyperbolicity:}

Certainly,  $DF({\bf w};H_b)$ is invertible  for any ${\bf w}\in\Sigma (\mathcal{A}_b,H_b)$ and $|b|<\epsilon_0$.
Suppose ${\bf w}=\theta(b;{\bf w}^\dag)$, then Theorem \ref{la23}(i) says that $\|{\bf w}-{\bf w}^\dag\|_\infty$ is uniformly bounded in ${\bf w}^\dag$. Therefore, the difference between the norms of $DF({\bf w}; H_b)$ and $DF({\bf w}^\dag; H_0)$ is also uniformly bounded in ${\bf w}^\dag$. This implies that $\|DF({\bf w}; H_b)^{-1}\|_\infty$ has a uniform bound in $\mathbf{w}$ by the uniform boundedness of $\|DF({\bf w}^\dag; H_0)^{-1}\|_\infty$ in ${\bf w}^\dag$. Now, the hyperbolicity follows from Theorem \ref{hypequiv}.   
 
The proof of Theorem \ref{mainthm} is complete. 
\qed  \\

\subsection{Proof of Theorem \ref{thm:RK278}}

Now, let $b_0$ satisfy  $|b_0|<\epsilon_0$, then the function 
\[
\Theta:\{b|~0\not= |b|<\epsilon_0\}\times \mathcal{A}_{b_0}\to\C^2, \qquad (b, \theta(b_0; {\bf w}^\dag)_0)\mapsto \theta(b;{\bf w}^\dag)_0,
\]
defines a holomorphic motion, where $\theta(\cdot; {\bf w}^\dag)$ is the holomorphic function given in \eqref{thetabw}. When $b\to 0$, the map $\Phi_b$ defined in \eqref{Phib} becomes an identity map. Because $\mathcal{A}$ is invariant under the restriction $H_0|_{\mathcal{A}}$ as shown in Proposition \ref{AH&Lambda}(i), we get 
\[
\mathcal{A}_b=\bigcup_{{\bf w}^\dag\in\Sigma(\mathcal{A}, H_0)}\pi_0\circ\Phi_b({\bf w}^\dag)\to
\bigcup_{{\bf w}^\dag\in\Sigma(\mathcal{A}, H_0)}\pi_0({\bf w}^\dag) = \bigcup_{w\in\mathcal{A}}w=\mathcal{A}
\]
in the Hausdorff topology as $b\to 0$.
\qed \\

\subsection{Proof of Theorem \ref{mainthm2}}

Let $(x^{b_0},y^{b_0})\in\mathcal{B}_{b_0}$, 
$(x^b,y^b):=\Theta(b, (x^{b_0}, y^{b_0}))\in\mathcal{B}_b$, and $\mathbf{z}^b=(z_i^b)_{i\in\Z}=(x_i^b,y_i^b)_{i\in\Z}$ with $(x_i^b,y_i^b)=H_b^i(x^b,y^b)$. It is clear that $(x_i^b, y_i^b)=\Theta(b, (x_i^{b_0},y_i^{b_0}))$ for all integer $i$.
To prove the theorem, we shall show that $\mathbf{z}^b\in\Sigma (\mathcal{B}_b, H_b)$ converges in the uniform topology to a point, say $\mathbf{z}^0\in \Sigma (\mathcal{B}, H_0)$, as $b\to 0$ with $\mathcal{B}$ having the desired property, and that $\mathbf{z}^{b_0}\mapsto \mathbf{z}^0$ is injective.

It has been known  that $J(H_b)$ is hyperbolic when $|b|$ is small \cite{BS1992, FS1992, HO-V1995,  IS2010}. Consequently, the linear operator $DF(\mathbf{z}^b; H_b)$ is invertible.
Indeed, for any neighbourhood $U$ of $\mathcal{A}=g(J(q))$, 
 we can take $|b|$ small enough so that $J(H_b)\subset U$. As a matter of fact, there is open set $\omega$ containing $J(q)$ and  $\delta>0$ such that $ \{(x,y)|~ y\in\omega~\mbox{and}~ |x-q(y)|<\delta\}$ is contained in $U$, and contains the Julia set $J(H_b)$ provided that $|b|$ is sufficiently small.
 For  $z=(x,y) \in\Gamma:=g(\C)$, the eigenvalue $2x$ of $DH_0(x,y)$ associates an eigenvector $(2x,1)$, while the eigenvalue $0$ associates eigenvector $(0,1)$. Thus,  the unstable cone $\mathcal{C}^u_z$ for a point $z\in J(H_b)$ can be a small conical neighbourhood of $T_z\mathcal{A}$ and can be extended continuously to  $U$, while  the stable cone $\mathcal{C}^s_z$ can be any conical neighbourhood of the vector $(0,1)$. Since $q$ is hyperbolic, there exists an adapted metric on $\Gamma$ such that $H_0$ is expanding on $\mathcal{A}$ along the tangent bundle $T\mathcal{A}$. We may extend this metric to $T\C^2|_\mathcal{A}$ by defining it in an arbitrary way on the normal bundle to $\Gamma$, and then extend this Riemannian metric continuously to $U$ in $\C^2$. Since we extended an adapted metric, and since $H_0$ has ``infinite contraction" along the direction $(0,1)$, the cone criterion holds on a neighbourhood of $J(H_b)$ when $|b|$ is small.
We refer the reader to \cite{BS1992} for more details. 

The expansion and contraction rates in the unstable and stable cones can be taken independent of $b$ as long as $|b|<\epsilon_2$ for some $0<\epsilon_2$. In view of Remarks \ref{InverseNorm} and \ref{InverseNormL=0}, there is $N>0$ such that $\|DF(\mathbf{z}^b; H_b)^{-1}\|_\infty<N$ for all $\mathbf{z}^b\in \Sigma(J(H_b), H_b)$ and $0<|b|<\epsilon_2$. By the implicit function theorem, $\mathbf{z}^b$ regarded as a function of $b$ satisfies
\[ \frac{d\mathbf{z}^b}{db}=-DF(\mathbf{z}^b; H_b)^{-1}\boldsymbol{\eta}
\]
 with $\boldsymbol{\eta}=(\eta_i)_{i\in\Z}$ the partial derivative $\partial_b F(\cdot; H_b)$  with respect to the parameter $b$ evaluated at $\mathbf{z}^b$. Direct calculation shows $\eta_i=(-y_i^b,0)$ for every integer $i$. It is known (e.g. \cite{Mumm2008}) that 
 $\sup_{(x, y)\in J(H_{b})}\max\left\{|x|,|y|\right\}\le (1+|b|+\sqrt{(1+|b|)^2+4|c|})/2=:R$.
This yields the estimate $ \|d\mathbf{z}^b/db\|_\infty\le\|DF(\mathbf{z}^b; H_b)^{-1}\|_\infty\cdot\|\boldsymbol{\eta}\|_\infty<NR$, and then leads to
\[ \|\mathbf{z}^{b_0}- \mathbf{z}^b\|<\begin{cases}
                                                         (b_0-b)NR & \mbox{for}~ 0<b<b_0<\epsilon_2 \\
                                  (b-b_0)NR & \mbox{for}~ -\epsilon_2<b_0<b<0.
\end{cases}
\]
Hence, the improper integral 
\[ 
    \mathbf{z}^0=\mathbf{z}^0(\mathbf{z}^{b_0}):=\lim_{b\to 0}\mathbf{z}^b=\mathbf{z}^{b_0}+\lim_{\tilde{b}\to 0}\int_{b_0}^{\tilde{b}} \frac{d\mathbf{z}^b}{db} db
\]
 exists. Because $F(\cdot; H_b)$ depends continuously on its variable and because $H_b(\cdot)$ depends continuously on its parameter $b$, we get $F(\mathbf{z}^0;H_0)=F(\lim_{b\to 0}\mathbf{z}^b; \lim_{b\to 0}H_b)=\lim_{b\to 0} F(\mathbf{z}^b; H_b)=0$.

It is well known that if $T: l_\infty(\Z, \C^2)\to  l_\infty(\Z, \C^2)$ is a linear operator such that  $\|T-DF(\mathbf{z}^b; H_b)\|_\infty<\|DF(\mathbf{z}^b; H_b)^{-1}\|_\infty^{-1}$ then $T$ is also invertible (see \cite{KF1970} for example). Now, 
\begin{eqnarray*}
 \|DF(\mathbf{z}^0; H_0)-DF(\mathbf{z}^{b_0}; H_{b_0})\|_\infty
  &=&\sup_{i\in\Z}\sup_{|(\xi_i, \tilde{\xi}_i)|=1}\left|\left((2x_i^0-2x_i^{b_0})\xi_i+b_0\tilde{\xi}_i, 0 \right)\right| \\
  &\le&\sup_{i\in\Z}2|x_i^0-x_i^{b_0}|+|b_0| \\
&<& |b_0|(2NR+1) \\
&<& 1/N
\end{eqnarray*}
provided that $|b_0|<\epsilon_1$ for some $\epsilon_1\le \min\left\{\epsilon_2, ~ N^{-1}(2RN+1)^{-1}\right\}$. Therefore, $DF(\mathbf{z}^0; H_0)$ is invertible. As a consequence of the implicit function theorem, $\mathbf{z}^{b_0}\mapsto \mathbf{z}^0$ is bijective. By setting $\epsilon_1$ smaller if needed and invoking   Theorem \ref{la23}(ii),  the mapping $\mathbf{z}^{b_0}\mapsto \mathbf{z}^0$ is a homeomorphism with the product topology.

Let $\mathcal{B}:=\bigcup_{\mathbf{z}^{b_0}\in\Sigma (\mathcal{B}_{b_0}, H_{b_0})}z^0_0(\mathbf{z}^{b_0})$. Then, $\mathcal{B}\subset \Gamma$ and is a compact invariant set of $H_0|_\mathcal{B}$. Let $\Lambda^\prime=g^{-1}(\mathcal{B})$, then $\Lambda^\prime$ is a compact  invariant  set of $q$. Since the unstable dimension of $\Lambda^\prime$ is one, $\Lambda^\prime$ is expanding and  must locate inside $J(q)$. \qed

\appendix
 
\section{Appendix}

In this appendix, we give a general setting which is the idea behind the proof of Theorem \ref{mainthm}.

Given a complex Banach space $l_\infty(\Z, \C^k)$, $k\ge 1$, denote by $\sigma$ the usual shift operator on it. 
Let $V$ be an open set in $l_\infty (\Z, \C^k)$,  $\mathcal{E}$ an open set in $\C^s$, $s\ge 1$, with the Euclidean metric, and $G:V\times\mathcal{E}\to l_\infty (\Z, \C^k)$ a holomorphic function. Let $\Sigma\subset V$ be a uniformly discrete subset of $l_\infty(\Z, \C^k)$ and $e^\dag$ a point in $\mathcal{E}$  such that $G({\bf v}^\dag, e^\dag)=0$ for every ${\bf v}^\dag\in\Sigma$. Denote by $B({\bf v}^\dag,\lambda)$ and  $\bar{B}({\bf v}^\dag,\lambda)$, respectively, the open and closed ball of radius $\lambda$ centered at ${\bf v}^\dag$ in $l_\infty (\Z,\C^k)$. Assume that the continuous linear operator $D_{\bf v}G({\bf v}^\dag,e^\dag):l_\infty (\Z,\C^k)\to l_\infty (\Z,\C^k)$, which is the partial derivative of $G$ at $({\bf v}^\dag, e^\dag)$ with respect to the first variable, is invertible and its inverse is likewise continuous for every ${\bf v}^\dag\in\Sigma$. By the implicit function theorem (e.g. \cite{Whit1965}), for every ${\bf v}^\dag$ there exists  $\hat{\delta}({\bf v}^\dag)$ and a unique holomorphic function $\theta(\cdot;{\bf v}^\dag):\mathcal{E}\to V$ such that $G(\theta(e ;{\bf v}^\dag),e)=0$ and $\theta(e^\dag;{\bf v}^\dag)={\bf v}^\dag$ provided $0\le |e-e^\dag|<\hat{\delta}({\bf v}^\dag)$.  Moreover, there exists $\hat{\lambda}({\bf v}^\dag)>0$ such that if $G({\bf v},e)=0$, $0\le |e-e^\dag|<\hat{\delta} ({\bf v}^\dag)$, 
and ${\bf v}\in \bar{B}({\bf v}^\dag,\hat{\lambda}({\bf v}^\dag))$, then ${\bf v}=\theta(e;{\bf v}^\dag)$. In \cite{CCY2014}, a theorem which  essentially is the same as  Theorem \ref{la23} below is developed for the theory of anti-integrability. (See also a version of  the theorem  for non-autonomous maps presented in \cite{Chen2005}.)  The theorem presented  in \cite{CCY2014} does not require $\Sigma$ to be uniformly discrete because it is automatically satisfied for $\Sigma$ being the set of anti-integrable orbits. However, the proof would be almost unchanged even if the uniform discreteness is considered.   For completeness sake, we include a proof of Theorem \ref{la23}.

\begin{theorem} \label{la23}
~
\begin{enumerate}
\item Assume    $\| D_{\bf v} G( {\bf v}^\dag, e^\dag)^{-1}\|_\infty$  is bounded above on $\Sigma\times\{e^\dag\}$, and assume for any $\gamma>0$ there exist $\lambda_0>0$ and $\delta_0>0$ such that for all ${\bf v}^\dag\in\Sigma$ we have $\|G({\bf v},e)\|_\infty<\gamma$ and $\|D_{\bf v}G({\bf v},e)-D_{\bf v}G({\bf v}^\dag,e^\dag)\|_\infty<\gamma$ whenever ${\bf v}\in B( {\bf v}^\dag, \lambda_0)$ and $| e- e^\dag|<\delta_0$. Then, there are $\delta_1$, $\lambda_1>0$ such that $\inf_{{\bf v}^\dag\in\Sigma}\hat{\delta}({\bf v}^\dag)\ge\delta_1$, $\inf_{{\bf v}^\dag\in\Sigma}\hat{\lambda}({\bf v}^\dag)\ge\lambda_1$, and the map 
\begin{eqnarray*}
        \Phi_e:\Sigma &\to&  \mathcal{I}:=\bigcup_{{\bf v}^\dag\in\Sigma}\theta(e; {\bf v}^\dag), \\
                           {\bf v}^\dag  &\mapsto& \theta(e;{\bf v}^\dag)\subset \bar{B}({\bf v}^\dag, \lambda_1),
\end{eqnarray*}
is a bijection provided $0\le |e-e^\dag|<\delta_1$.  
\item  There exists $0<\delta_2 \le \delta_1$ such that with the product topology the map $\Phi_e$ is a continuous function  thus  a homeomorphism from $\Sigma$ to $\mathcal{I}$ provided that $0\le |e-e^\dag|< \delta_2$ and that $\Sigma$ is compact (within the product topology).
\item In addition to the assumptions in (i), if $\Sigma$ is forward $\sigma$-invariant, i.e. $\sigma(\Sigma)=\Sigma$, and if $G(\cdot,e)$ commutes with $\sigma$ for all $e\in\mathcal{E}$, i.e. $\sigma\circ G({\bf v},e)=G(\sigma({\bf v}),e)$, then the following diagram commutes
 \[
  \begin{matrix}
        \Sigma & \mapright{\sigma} & \Sigma \cr
          \mapldown{\Phi_e} & &  \maprdown{\Phi_e} \cr
    {\mathcal{I}} & \mapright{\sigma} & {\mathcal{I}}
\end{matrix}
\]
provided $0\le |e-e^\dag|<\delta_1$.  
\end{enumerate}
\end{theorem}
\proof
 (i) For each ${\bf v}^\dag\in\Sigma$ and $e\in\mathcal{E}$, define a map $\mathcal{G}(\cdot;{\bf v}^\dag,e):V\to l_\infty(\Z, \C^k)$, ${\bf v}\mapsto {\bf v}-D_{\bf v}G({\bf v}^\dag,e^\dag)^{-1}G({\bf v},e)$. By assumptions in assertion (i), there are $0<\lambda_1$ and $0<\delta_1$ such that for ${\bf u}$ in closed ball $\bar{B}({\bf v}^\dag, \lambda_1)$ and $|e-e^\dag|\le\delta_1$, we have
\[ \|D_{\bf v}G({\bf v}^\dag,e^\dag)^{-1}\|_\infty ~\|D_{\bf v}G({\bf u},e)-D_{\bf v}G({\bf v}^\dag,e^\dag)\|_\infty \le 1/2\]
and
\begin{equation}
 \|D_{\bf v}G({\bf v}^\dag,e^\dag)^{-1}\|_\infty ~\|G({\bf v}^\dag,e)\|_\infty <\lambda_1 /2. \label{etadelta}
\end{equation}
Thus, for ${\bf v}$, ${\bf w}\in \bar{B}({\bf v}^\dag,\lambda_1)$ and $|e-e^\dag|\le\delta_1$, we get 
\begin{eqnarray}
 \lefteqn{ \left\|(\mathcal{G}({\bf v};{\bf v}^\dag,e)-\mathcal{G}({\bf w};{\bf v}^\dag,e)\right\|_\infty  } \nonumber \\
&\le &  \sup_{\bar{\bf u}\in \{{\bf w}+t({\bf v}-{\bf w}):~ 0\le t\le 1\}}\left\|D_{\bf v}G({\bf v}^\dag,e^\dag)^{-1} \left(
                          D_{\bf v}G({\bf v}^\dag,e^\dag)-  D_{\bf v}G(\bar{\bf u},e)  
                                                                                 \right)\right\|_\infty\cdot\|({\bf v}-{\bf w})\|_\infty\nonumber\\
&\le& \|{\bf v}-{\bf w}\|_\infty /2 \label{401}
\end{eqnarray}
and then
\begin{eqnarray}
 \left\|\mathcal{G}({\bf v};{\bf v}^\dag,e)-{\bf v}^\dag\right\|_\infty
&\le & \left\| \mathcal{G}({\bf v};{\bf v}^\dag,e)-\mathcal{G}({\bf v}^\dag;{\bf v}^\dag,e)\right\|_\infty+\left\|\mathcal{G}({\bf v}^\dag;{\bf v}^\dag,e)-{\bf v}^\dag\right\|_\infty\nonumber\\
&\le& \|{\bf v}-{\bf v}^\dag\|_\infty /2+\|D_{\bf v}G({\bf v}^\dag,e^\dag)^{-1}\|_\infty ~\|G({\bf v}^\dag,e)\|_\infty\label{402}\\
&<&\lambda_1. \nonumber
\end{eqnarray}
This implies that $\mathcal{G}(\cdot;{\bf v}^\dag,e)$ is a contraction map with contraction constant at least $1/2$ on $\bar{B}({\bf v}^\dag,\lambda_1)$ for any ${\bf v}^\dag\in\Sigma$
 and $|e-e^\dag|\le\delta_1$. Hence, $\inf_{{\bf v}^\dag\in\Sigma}\hat{\delta}({\bf v}^\dag)\ge\delta_1>0$.

The radius $\lambda_1$ is independent of ${\bf v}^\dag$ and $e$, and $\Phi_e ({\bf v}^\dag)$ is the unique fixed point in $\bar{B}({\bf v}^\dag,\lambda_1)$ for $\mathcal{G}(\cdot;{\bf v}^\dag,e)$.  Hence $\inf_{{\bf v}^\dag\in\Sigma}\hat{\lambda}({\bf v}^\dag)\ge\lambda_1>0$. Because $\Sigma$ is uniformly discrete, the balls $\bar{B}({\bf v}^\dag,\lambda_1)$, ${\bf v}^\dag\in\Sigma$, are disjoint in $l_\infty (\Z, \C^k)$ provided that $\lambda_1$ is sufficiently small. It follows that $\Phi_e$ is bijective on $\Sigma$.

(ii) There exists $0<\delta_2\le\delta_1$ such that $\|D_{\bf v}G({\bf v}^\dag,e^\dag)^{-1}\|_\infty ~ \|G({\bf v}^\dag,e)\|_\infty<\lambda_1/4$ provided $|e- e^\dag|\le \delta_2$. Thence, in view of \eqref{401}-\eqref{402}, the ball $\bar{B}({\bf v}^\dag,\lambda_1/2)$ is mapped into itself under $\mathcal{G}(\cdot; {\bf v}^\dag,e)$, so $\|\Phi_e ({\bf v}^\dag)-{\bf v}^\dag\|_\infty<\lambda_1/2$ for all ${\bf v}^\dag\in\Sigma$ and $0\le |e-e^\dag|<\delta_2$. Given ${\bf v}^\dag_*\in \Sigma$ and  $\left({\bf v}^{\dag (k)}\right)_{k\ge 0}$  a convergent sequence to ${\bf v}^\dag_*$ in $\Sigma$ (with  the product topology),  suppose that $\left(\Phi_e({\bf v}^{\dag (k)})\right)_{k\ge 0}$ converges, via a subsequence if necessary, to a point ${\bf v}^*$ in
 $\bigcup_{{\bf v}^\dag\in\Sigma}\bar{B}({\bf v}^\dag, \lambda_1/2)$. (Note that  $\bigcup_{{\bf v}^\dag\in\Sigma}\bar{B}({\bf v}^\dag, \lambda_1/2)$ is compact in the product topology.) Assume the subsequence. For any $N\in\N$, there is $K\in\N$ such that 
\[ |v_i^{\dag (k)}-v_i^\dag|<\lambda_1/2\]
for all $k>K$,  $|i|<N$, 
 and then 
\begin{eqnarray*}
|\Phi_e({\bf v}^{\dag (k)})_i-v_i^\dag| &\le& |\Phi_e({\bf v}^{\dag (k)})_i-v_i^{\dag (k)}|+|v_i^{\dag (k)}-v_i^\dag|\\
                                                                    &<&\lambda_1
\end{eqnarray*}
for $0\le |e-e^\dag|<\delta_2$. Passing to $N\to\infty$, we have 
\[
  \|{\bf v}^* -{\bf v}^\dag_*\|_\infty\le\lambda_1.
\]
Besides, 
\[ G({\bf v}^*,e)=G(\lim_{k\to\infty}\Phi_e({\bf v}^{\dag (k)}),e)=\lim_{k\to\infty}G(\Phi_e({\bf v}^{\dag (k)}),e)=0\]
because $G(\cdot,e)$ is continuous on $\bigcup_{{\bf v}^\dag\in\Sigma}\bar{B}({\bf v}^\dag, \lambda_1)$, which is contained in $V$. This means that ${\bf v}^*$ is a zero of $G(\cdot,e)$ in $\bar{B}({\bf v}^\dag_*,\lambda_1)$. From the assertion (i), we conclude that ${\bf v}^*$ must be $\Phi_e({\bf v}_*^{\dag})$.

(iii) $G(\cdot,e)$ has  a unique zero at $\Phi_e({\bf v}^\dag)$ in $B({\bf v}^\dag,\lambda_1)$, so does  $G(\sigma(\cdot),e)$
 since $\sigma\circ G(\cdot,e)=G(\sigma(\cdot),e)$. This implies that $G(\cdot,e)$ has a unique zero at $\sigma\circ\Phi_e({\bf v}^\dag)$ in $B(\sigma({\bf v}^\dag),\lambda_1)$. Because $G(\cdot,e)$ has been shown to have a unique zero at $\Phi_e (\sigma({\bf v}^\dag))$ in $B(\sigma({\bf v}^\dag),\lambda_1)$, it must have $\sigma\circ\Phi_e ({\bf v}^\dag)=\Phi_e(\sigma({\bf v}^\dag))$ by the uniqueness.
\qed

\section*{Acknowledgments}

This work was partly supported by NSC 101-2115-M-001-010, MOST 107-2115-M-001-006 and 109-2115-M-001-006.
 The author is grateful to  Zin Arai, Yutaka Ishii and Tomoki Kawahira  for useful conversations when an early version of this paper  was announced in the RIMS workshop ``New Developments in Complex Dynamical Systems" (Kyoto, 2012) and in the ``Complex Dynamics Seminar" at Nagoya University (Nagoya, 2012). He also thanks Kenneth Palmer for discussing on exponential dichotomy, and thanks for  the hospitality of Kyoto, Kyushu, Nagoya, and Suzhow Universities during his visits.

\end{document}